\documentclass[12pt,a4paper,english,intoc,bibliography=totoc,index=totoc,BCOR10mm,captions=tableheading,titlepage,fleqn]{scrbook}
\usepackage{lmodern}

\usepackage[T1]{fontenc}
\usepackage[latin9]{inputenc}
\usepackage{babel}
\usepackage{amsmath}
\usepackage{amssymb}
\usepackage[unicode=true,
 bookmarks=true,bookmarksnumbered=true,bookmarksopen=true,bookmarksopenlevel=2,
 breaklinks=false,pdfborder={0 0 0},backref=false,colorlinks=false]
 {hyperref}
\hypersetup{pdftitle={Abstract},
 pdfauthor={Uwe Stöhr},
 pdfpagelayout=OneColumn, pdfnewwindow=true, pdfstartview=XYZ, plainpages=false}
\usepackage{breakurl}

\makeatletter

\usepackage{calc}

\makeatother

\begin{document}
\noindent \begin{center}
\textbf{\Large Madness and regularity properties}
\par\end{center}{\Large \par}

\noindent \begin{center}
{\large Haim Horowitz and Saharon Shelah}
\par\end{center}{\large \par}

\noindent \begin{center}
\textbf{Abstract}
\par\end{center}

\noindent \begin{center}
Starting from an inaccessible cardinal, we construct a model of $ZF+DC$
where there exists a mad family and all sets of reals are $\mathbb Q$-measurable
for $\omega^{\omega}$-bounding sufficiently absolute forcing notions
$\mathbb Q$.%
\footnote{Date: May 16, 2017

2000 Mathematics Subject Classification: 03E35, 03E15, 03E25

Keywords: mad families, Lebesgue measurability, amalgamation

Publication 1113 of the second author

Partially supported by European Research Council grant 338821 and
NSF grant DMS-1362974.%
}
\par\end{center}

\textbf{\large Introduction}{\large \par}

Our study concerns the interactions between mad families and other
types of pathological sets of reals. Given a forcing notion $\mathbb Q$
whose conditions are subtrees of $\omega^{<\omega}$ ordered by reverse
inclusion, the notion of $\mathbb Q-$measurability is naturally defined.
As the existence of mad families and non-$\mathbb Q-$measurable sets
follows from the axiom of choice, one may consider the possible implications
between the existence of mad families and the existence of non-$\mathbb Q-$measurable
sets. The study of models of $ZF+DC$ where no mad families exist
was initiated by Mathias in {[}Ma{]}, more results were obtained recently
in {[}HwSh1090{]}, {[}NN{]} and {[}To{]}. Models of $ZF+DC$ where
all sets of reals are $\mathbb Q$-measurable for various forcing
notions $\mathbb Q$ were first studied by Solovay in {[}So{]}.

Our main goal is to show that $\mathbb Q-$measurability for $\omega^{\omega}$-boundning
sufficiently absolute forcing notions does not imply the non-exsitence
of mad families. In particular, as Random real forcing is $\omega^{\omega}$-bounding,
it will follow that Lebesgue measurability for all sets of reals does
not imply the non-existence of mad families.

We follow the strategy of {[}Sh218{]}, where a model of $ZF+DC+"$all
sets of reals are Lebesgue measurable but there is a set without the
Baire Property$"$ was constructed. Fixing an inaccessible cardinal
$\kappa$, we define a partial order $AP$ consisting of pairs $(\mathbb P,\Gamma)$,
where $\mathbb P$ is a forcing notion from $H(\kappa)$ and $\Gamma$
is an approximation of the desired mad family such that finite unions
of members of $\Gamma$ are not dominated by reals from $V$. We shall
obtain our model by forcing with this partial order and then with
the partial order introduced generically by $AP$. The main point
will be an amalgamation argument for $AP$ (over $\mathbb Q$-generic
reals for an appropriate $\mathbb Q$), which will allow us to repeat
Solovay's argument from {[}So{]}.

Remark: It was brought to our attention by Paul Larson and Jindra
Zapletal that a model of {}``every set of reals is Lebesgue measurable
and there is a mad family'' can also be constructed using the arguments
from Section 5 of their paper {[}LZ{]}. However, they assume the existence
of a proper class of Woodin cardinals, while in this paper we only
assume the existence of an inaccessible cardinal.

$\\$

\textbf{\large The main result}{\large \par}

\textbf{Hypothesis 1: }Throughout the paper, $\bold f$ will be a
fixed forcing frame (defined below) with $\kappa_{\bold f}=\kappa$
a fixed inaccessible cardinal.

\textbf{Definition 2: }Let $\bold f=(\kappa_{\bold f},\bold{P}_{\bold f},\bold{Q}_{\bold f})=(\kappa,\bold P,\bold Q)$
be a forcing frame when:

a. $\kappa$ is the inaccessible cardinal from Hypothesis 1.

b. $\bold P$ is the set of forcing notions from $H(\kappa)$.

c. $\bold Q$ is a family of $\omega^{\omega}-$bounding forcing notions
with sufficiently absolute definitions.

d. If $\mathbb P \in \bold P$ and $V^{\mathbb P} \models "\mathbb Q \in \bold Q"$,
then $\mathbb Q \in H(\kappa)^{(V^{\mathbb P})}$.

\textbf{Definition 3: }Let $AP=AP_{\bold f}$ be the partial order
defined as follows:

a. $a\in AP$ iff $a$ has the form $(\mathbb P,\Gamma)=(\mathbb{P}_a,\Gamma_a)$
where:

1. $\mathbb P \in \bold P$ and $\Gamma$ is an infinite set of canonical
$\mathbb P-$names of reals such that $\Vdash_{\mathbb P} "\Gamma$
is almost disjoint$"$.

2. If $\underset{\sim}{\tau} \in \Gamma$, then $\Vdash_{\mathbb P} "\underset{\sim}{\tau}$
is an infinite subset of $\omega"$.

3. For $a\in AP$, let $\Omega_a$ be the set of $\underset{\sim}{\tau} \in \Gamma_a$
which are objects and not just names.

4. If $1\leq n$, $\underset{\sim}{a_0},...,\underset{\sim}{a_{n-1}} \in \Gamma_a \setminus \Omega_a$,
$\underset{\sim}{a}=\underset{l<n}{\cup}\underset{\sim}{a_l}$ and
$f_{\underset{\sim}{a}} : \omega \rightarrow \omega$ is the function
enumerating $\underset{\sim}{a}$ in an increasing order, then $\Vdash_{\mathbb P} "f_{\underset{\sim}{a}}$
is not dominated by any $f\in (\omega^{\omega})^V"$.

b. $a\leq_{AP} b$ iff

1. $\mathbb{P}_a \lessdot \mathbb{P}_b$.

2. $\Gamma_a \subseteq \Gamma_b$.

3. If $\underset{\sim}{a_0},...,\underset{\sim}{a_{n-1}} \in \Gamma_b \setminus \Gamma_a$,
$\underset{\sim}{a}=\underset{l<n}{\cup} \underset{\sim}{a_l}$ and
$f_{\underset{\sim}{a}}$ enumerates $\underset{\sim}{a}$ in an increasing
order, then $\Vdash_{\mathbb{P}_b} "f_{\underset{\sim}{a}}$ is not
dominated by any member of $(\omega^{\omega})^{V[G\cap \mathbb{P}_a]}$.

\textbf{Observation 4: }$(AP,\leq)$ is indeed a partial order.

\textbf{Proof: }Suppose that $a\leq b$ and $b\leq c$. Let $\underset{\sim}{a_0},...,\underset{\sim}{a_{n-1}} \in \Gamma_c \setminus \Gamma_a$,
and let $\underset{\sim}{a}$ and $f_{\underset{\sim}{a}}$ be as
in Definition 3(b)(3). We may assume wlog that for some $0<m<n$,
$\underset{\sim}{a_l} \in \Gamma_b$ iff $l<m$ (the cases $m=0$
and $m=n$ are trivial). Let $G_c \subseteq \mathbb{P}_c$ be $V$-generic
and let $G_a=G_c \cap \mathbb{P}_a$ and $G_b=G_c \cap \mathbb{P}_b$.
Let $g=(n_i : i<\omega) \in V[G_a]$, wlog $g$ is increasing. We
shall prove that $f_{\underset{\sim}{a}}$ is not dominated by $g$.

Let $a_i=\underset{\sim}{a_i}[G_c]$, $a=\underset{\sim}{a}[G_c]$
and $b=\underset{l<m}{\cup}a_l$.

Subclaim 1: For infinitely many $i$, $[n_i,n_{i+1}) \cap (\underset{l<n}{\cup}a_l)=\emptyset$.

Subclaim 2: Subclaim 1 is equivalent to $"f_a$ is not dominated by
$g"$.

Proof of Subclaim 1: Let $u=\{i : [n_i,n_{i+1}) \cap b=\emptyset\} \in V[G_b]$.
By the fact that $a\leq b$ and by subclaim 2, $u$ is infinite. Let
$(i(l) : l<\omega) \in V[G_b]$ be an increasing enumeration of $u$,
so $(n_{i(l)} : l<\omega) \in V[G_b]$ is increasing. Let $c=\underset{m\leq l<n-1}{\cup}a_l$
and $v=\{l : [n_{i(l)},n_{i(l+1)}) \cap c=\emptyset\}$. As before,
$v$ is infinite. If $l \in v$ then $c\cap [n_{i(l)},n_{i(l+1)})=\emptyset$
and therefore, $c\cap [n_{i(l)},n_{i(l)+1})=\emptyset$. Similarly,
if $l\in v$ then $i(l) \in u$ and therefore $b\cap [n_{i(l)},n_{i(l)+1})$.
It follows that $l\in v \rightarrow (b\cup c) \cap [n_{i(l)},n_{i(l)+1})=\emptyset$,
and as $v$ is infinite, we're done.

Proof of Subclaim 2: Suppose that $f_{\underset{\sim}{a}}$ is not
dominated by any $g\in (\omega^{\omega})^{V^{\mathbb{P}_a}}$ and
let $g=(n_i : i<\omega) \in V^{\mathbb{P}_a}$ be increasing. Choose
$f\in V^{\mathbb{P}_c}$ such that $f$ is increasing, $l<f(l)$ for
every $l$ and $|\{i : n_i \in [l,f(l))\}|$ is sufficiently large
(e.g. $>2^l$). By our assumption, for infinitely many $l$, $f(l)\leq$
the $l$th member of $\underset{\sim}{a}$, and therefore $|\underset{\sim}{a} \cap f(l)| \leq l$.
Let $u=\{ l : |\underset{\sim}{a} \cap f(l)| \leq l\}$, so $u$ is
infinite. For $l\in  u$, $l+1<|\{i : l\leq i, [n_i,n_{i+1}) \subseteq [l,f(l))\}|$,
and as $u$ is infinite, for some $i$ such that $l\leq i$, $[n_i,n_{i+1}) \subseteq [l,f(l))$
and $[n_i,n_{i+1}) \cap \underset{\sim}{a}=\emptyset$. Therefore,
for infinitely many $i$, $[n_i,n_{i+1}) \cap \underset{\sim}{a}=\emptyset$.

In the other direction, suppose that $f_{\underset{\sim}{a}}$ satisfies
the condition of Subclaim 1. Let $g\in (\omega^{\omega})^{V^{\mathbb{P}_a}}$,
we shall prove that $f_{\underset{\sim}{a}}$ is not dominated by
$g$. We may assume wlog that $g$ is increasing. Choose the sequence
$(n_i : i<\omega)$ by induction such that $n_0=0$ and $n_{i+1}>n_i+g(n_i)$,
so $(n_i : i<\omega) \in V^{\mathbb{P}_a}$. By the assumption, the
set $u=\{i : [n_i,n_{i+1}) \cap \underset{l<n}{\cup}a_l=\emptyset\}$
is infinite. For every $i\in u$, $|a\cap n_i| \leq n_i$, therefore
$n_i < f_a(n_i)$. As $[n_i,n_{i+1}) \cap a=\emptyset$, it follows
that $n_{i+1} \leq f_a(n_i)$, therefore $g(n_i)<n_{i+1} \leq f_a(n_{n_i})$,
so $f_a$ is not dominated by $g$. $\square$

\textbf{Observation 4: }a. Every $\mathbb P \in \bold P$ is $\kappa-cc$,
and $\bold P$ is closed under $\lessdot-$increasing unions of length
$<\kappa$.

b. If $\mathbb P \in \bold P$ and $\underset{\sim}{\mathbb Q}$ is
a canonical $\mathbb{P}-$name of a case of $\bold Q$ which is in
$H(\kappa)$, then $\mathbb P \star \underset{\sim}{\mathbb Q} \in \bold P$.
$\square$

\textbf{Observation 5: }a. If $a\in AP$ then $(\{0\},\Omega_a) \in AP$
and $(\{0\},\Omega_a) \leq a$.

b. $AP$ is $(<\kappa)-$complete. $\square$

\textbf{Claim 6: }$(AP,\leq)$ has the division property, namely,
if $a\leq b$ and $\underset{\sim}{x}$ is a $\mathbb{P}_b$-name
of a real such that $\Vdash_{\mathbb{P}_b} "(\omega^{\omega})^{V[\mathbb{P}_a]}$
is cofinal in $(\omega^{\omega})^{V[\mathbb{P}_a,\underset{\sim}{x}]}"$,
then there is $a_1 \in AP$ such that:

a. $a\leq a_1 \leq b$.

b. $\Gamma_{a_1}=\Gamma_a$.

c. $\mathbb{P}_{a_1}=\mathbb{P}_a \star \underset{\sim}{x}$ in the
natural sense. $\square$

\textbf{Claim 7 ($(AP,\leq)$ has the amalgamation property): }Assume
that $a_0 \leq a_l$ $(l=1,2)$, then there are $b_l$ $(l\leq 3)$
and $g_l$ $(l\leq 2)$ such that:

a. $b_0 \leq b_l \leq b_3$ $(l=1,2)$.

b. $g_l$ is an isomorphism from $b_l$ to $a_l$.

c. $g_0 \subseteq g_l$ $(l=1,2)$. 

\textbf{Proof: }We may assume wlog that $\mathbb{P}_{a_0}$ is trivial
and that $\Omega_{a_1}=\Omega_{a_2}=\Gamma_{a_0}$ (as we can simply
take the quotients).

We define $\mathbb{P}_{b_3}$ as follows:

a. $p\in \mathbb{P}_{b_3}$ iff $p=(p_1,p_2) \in \mathbb{P}_{a_1} \times \mathbb{P}_{a_2}$
and for some $l(p), n_p, A_{p,1},A_{p,2},\underset{\sim}{a_{p,1}},\underset{\sim}{a_{p,2}}$
the following hold:

1. $l(p) \in \{1,2\}$ and $n_p <\omega$.

2. $A_{p,l}$ is a finite subset of $\Gamma_{a_l}$ with union $\underset{\sim}{a_{p,l}}$
$(l=1,2)$.

3. For every $n>n_p$, there is $r_{n} \in \mathbb{P}_{a_{l(p)}}$
such that $\mathbb{P}_{a_{l(p)}} \models p_{l(p)} \leq r_{n}$ and
$r_{n} \Vdash "\underset{\sim}{a_{p,l(p)}} \cap n \subseteq n_p"$.

b. $\mathbb{P}_{b_3} \models p\leq q$ iff

1. $p=(p_1,p_2), q=(q_1,q_2) \in \mathbb{P}_{a_1} \times \mathbb{P}_{a_2}$.

2. $p_l \leq q_l$ $(l=1,2)$.

3. $n_p \leq n_q$.

4. $A_{p,l} \subseteq A_{q,l}$ $(l=1,2)$.

5. There is no $n \in [n_p,n_q)$ such that $q_1 \nVdash "n\notin \underset{\sim}{a_{q,1}}"$
and $q_2 \nVdash "n\notin \underset{\sim}{a_{q,2}}"$.

We shall now define embeddings $f_l: \mathbb{P}_{a_l} \rightarrow \mathbb{P}_{b_3}$
$(l=1,2)$ as follows: For $p\in \mathbb{P}_{a_l}$, $f_l(p)=q \in \mathbb{P}_{b_3}$
will be the condition defined as follows:

a. $q_l=p$ and $q_{3-l}=0_{\mathbb{P}_{a_{3-l}}} \in \mathbb{P}_{a_{3-l}}$.

b. $l(q)=l$, $n_p=0$.

c. $A_{q,1}=\emptyset=A_{q,2}$.

Subclaim 0: $\mathbb{P}_{b_3}$ is a partial order.

Subclaim 1: For every $p=(p_1,p_2)\in \mathbb{P}_{b_3}$ and open
dense $I\subseteq \mathbb{P}_{3-l(p)}$, there is $q\in \mathbb{P}_{b_3}$
above $p$ such that $l(q)=3-l(p)$ and $q_{3-l(p)} \in I$.

Proof: Let $i=3-l(p)$ and let $p_i' \in I$ be above $p_i$. By the
definition of $AP$, $f_{\underset{\sim}{a_{p,i}}}$ is not dominated
by any function from $V$. We shall prove that there are $q_i \in \mathbb{P}_{a_i}$
above $p_i'$ and $n_*>n_p$ such that for every $n>n_*$, there is
$q'$ above $q_i$ such that $q'\Vdash "\underset{\sim}{a_{p,i}} \cap [n_*,n)=\emptyset"$.
Actually, $q_i=p_i'$ should work. Suppose not, then for every $n_*>n_p$
there is $n>n_*$ such that there is no $q'$ above $q_i$ forcing
that $\underset{\sim}{a_{p,i}} \cap [n_*,n)=\emptyset$. Now choose
$(n_j : j<\omega)$ by induction on $j$ as follows: $n_0=n_p+1$,
and $n_{j+1}$ is the minimal $n>n_j$ such that there is no $q'$
above $q_i$ forcing that $\underset{\sim}{a_{p,i}} \cap [n_j,n)=\emptyset$.
By the same argument as in the proof of observation 4, as $(n_j : j<\omega) \in V$,
$p_i' \Vdash "\underset{\sim}{a_{p,i}} \cap [n_j,n_{j+1})=\emptyset$
for infinitely many $j"$. Therefore, there is $q'$ above $p_i'$
and $i_*$ such that $q' \Vdash "\underset{\sim}{a} \cap [n_{i_*},n_{i_*+1})=\emptyset"$,
contradicting the choice of $n_{i_*+1}$. 

Now define $q\in \mathbb{P}_{b_3}$ as follows:

1. $q_i$ is as above.

2. $q_{l(p)}$ is any member of $\mathbb{P}_{l(p)}$ which is above
$p_{l(p)}$ and forces that $[n,n_*) \cap \underset{\sim}{a_{p,l(p)}}=\emptyset$
(such condition exists by clause (a)(3) in the definition of $\mathbb{P}_{b_3}$).

3. $l(q)=i$.

4. $n_q=n_*$.

5. $A_{q,l}=A_{p,l}$ and $\underset{\sim}{a_{q,l}}=\underset{\sim}{a_{p,l}}$
for $l=1,2$.

It's now easy to check that $q$ is as required.

Subclaim 2: a. $\{p \in \mathbb{P}_{b_3} : l(p)=i\}$ is dense in
$\mathbb{P}_{b_3}$ for $i=1,2$.

b. $I_n:=\{p \in \mathbb{P}_{b_3} : n_p>n\}$ is dense in $\mathbb{P}_{b_3}$.

Proof: (a) follows from Subclaim 1. (b) follows from the proof of
Subclaim 1, as we note that $n_q=n_*>n_p$ in that proof.

Subclaim 3: $f_l: \mathbb{P}_{a_l} \rightarrow \mathbb{P}_{b_3}$
is a complete embedding for $l=1,2$.

Proof: It suffices to show that $f_l$ is a complete embedding into
$\{p\in \mathbb{P}_{b_3} : l(p)=l\}$, which follows from the existence
of a projection $\pi: \{p\in \mathbb{P}_{b_3} : l(p)=l\} \rightarrow \mathbb{P}_{a_l}$
defined in the natural way.

Subclaim 4: For every finite $A_1 \subseteq \Gamma_{a_1}$ and $A_2 \subseteq \Gamma_{a_2}$,
the set $\{p \in \mathbb{P}_{b_3} : \underset{i=1,2}{\wedge}A_i \subseteq A_{p,i}\}$
is open dense.

Proof: In order to prove the claim by induction on $|A_1|+|A_2|$,
it suffices to prove it when $A_i=\{\underset{\sim}{b}\}$ and $A_{3-i}=\emptyset$
for $i\in \{1,2\}$. Let $p\in \mathbb{P}_{b_3}$ and suppose that
$l(p)=3-i$, it's now easy to extend $p$ simply by adding $\underset{\sim}{b}$
to $A_{p,i}$. If $l(p)=i$, then by previous claims, there is $q$
above $p$ such that $l(q)=3-l(p)$, and now extend $q$ as in the
previous case.

Subclaim 5: Let $\Gamma:=f_1(\Gamma_{a_1}) \cup f_2(\Gamma_{a_2})$,
then $\Gamma$ is a set of canonical $\mathbb{P}_{b_3}$-names of
infinite subsets of $\omega$ and $\Vdash_{\mathbb P} "\Gamma$ is
almost disjoint$"$.

Proof: The first part follows by the fact that $f_1$ and $f_2$ are
complete embeddings. In order to prove the second part, it suffices
to show that if $\underset{\sim}{r} \in \Gamma_{a_1}$ and $\underset{\sim}{s} \in \Gamma_{a_2}$,
then $\Vdash_{\mathbb P} "|\underset{\sim}{r} \cap \underset{\sim}{s}|< \aleph_0"$.
Given $p \in \mathbb{P}_{b_3}$, by Subclaim 4, there is a stronger
condition $q$ such that $\underset{\sim}{r} \in A_{q,1}$ and $\underset{\sim}{s} \in A_{q,2}$.
We shall prove that $q\Vdash "|\underset{\sim}{r} \cap \underset{\sim}{s}|<\aleph_0"$.

Recall that for every $n$, the set $I_n=\{r \in \mathbb{P}_{b_3} : n\leq n_r\}$
is dense. Now let $G\subseteq \mathbb{P}_{b_3}$ be generic over $V$
such that $q\in G$, then for every $n_q <n$, there is $q_n \in G$
such that $n\leq n_{q_n}$. By the definition of the partial order
$\leq_{\mathbb{P}_{b_3}}$ (clause (b)(5)), it follows that $q\Vdash_{\mathbb{P}_{b_3}} "|\underset{\sim}{r} \cap \underset{\sim}{s}|<\aleph_0"$. 

Subclaim 6: Let $b_3=(\mathbb{P}_{b_3},\Gamma_{b_3})$ where $\Gamma_{b_3}$
is $\Gamma$ from the previous subclaim, then $b_3$ satisfies clauses
$(1)+(2)$ from Definition (3)(a). As $\Omega_{a_1}=\Omega_{a_2}$,
it follows that $\Omega_{b_3}=\Omega_{a_1}=\Omega_{a_2}$.

For $l=1,2$, let $b_l=f_l(a_l) \in AP$, then clauses $(1)+(2)$
from Definition (3)(b) hold for $b_l$ and $b_3$. 

Subclaim 7: $b_3 \in AP$.

Proof: Let $A\subseteq \Gamma_{b_3} \setminus \Omega_{b_3}$ be finite,
so there are finite sets $A_l \subseteq \Gamma_{a_l} \setminus \Omega_{a_l}$
$(l=1,2)$ such that $A=f_1(A_1) \cup f_2(A_2)$. Let $(n_i : i<\omega) \in (\omega^{\omega})^V$
be increasing and let $\underset{\sim}{u}=\{i : [n_i,n_{i+1}) \cap (\cup \{ \underset{\sim}{a} : \underset{\sim}{a} \in A \})=\emptyset\}$.
Let $(p_1,p_2) \in \mathbb{P}_{b_3}$ and $n<\omega$, we shall find
$(q_1,q_2)$ and $i>n$ such that $(p_1,p_2) \leq (q_1,q_2) \in \mathbb{P}_{b_3}$
and $(q_1,q_2) \Vdash_{\mathbb{P}_{b_3}} "i\in \underset{\sim}{u}"$.
Without loss of generality, $l((p_1,p_2))=2$, and by Subclaim 4,
wlog $A_i \subseteq A_{(p_1,p_2),i}$ $(i=1,2)$. For $l=1,2$, let
$\underset{\sim}{a_l}=\cup \{ \underset{\sim}{a} : \underset{\sim}{a} \in A_l\}$,
so $\underset{\sim}{a_l}$ is a $\mathbb{P}_{a_l}-$name and $\Vdash_{\mathbb{P}_{a_l}}"(\exists^{\infty}i)(\underset{\sim}{a_l} \cap [n_i,n_{i+1})=\emptyset)"$.
Choose $(p_{1,l},j_{1,l} : l<\omega)$ by induction on $l<\omega$
such that:

1. $p_{1,0}=p_1$.

2. $\mathbb{P}_{a_1} \models p_{1,l} \leq p_{1,l+1}$.

3. $j_{1,l}>l+\underset{k<l}{\Sigma}j_{1,k}$.

4. $p_{1,l+1} \Vdash_{\mathbb{P}_{a_1}} "\underset{\sim}{a_1} \cap [n_{j_{1,l}},n_{j_{1,l}+1})=\emptyset"$.

For $l<\omega$, let $m_l=n_{j_{1,l}}$, so $(m_l : l<\omega) \in (\omega^{\omega})^V$
is increasing. Let $j$ be the minimal $j>n$ such that $n_{(p_1,p_2)} \leq m_j$.
By the proof of Subclaim 1, there are $p_1'$ above $p_{1,j+1}$ and
$k_*>n_{(p_1,p_2)}$ such that for every $k>k^*$ there is $p''$
above $p_1'$ forcing that $\underset{\sim}{a_{(p_1,p_2),1}} \cap [k^*,k)=\emptyset$.
As $l((p_1,p_2))=2$, there is $p_2'$ above $p_2$ forcing that $\underset{\sim}{a_{(p_1,p_2),2}} \cap [n_{(p_1,p_2)},k^*+m_{j+1})=\emptyset$.
Now let $(q_1,q_2)=(p_1',p_2')$, $n_{(q_1,q_2)}=k^*$, $l((q_1,q_2))=1$,
$A_{(q_1,q_1),i}=A_{(p_1,p_2),i}$ $(i=1,2)$, it's easy to see that
$(q_1,q_2)$ and $j$ are as required.

Subclaim 8: $b_l \leq b_3$ where $b_l=f_l(a_l)$ $(l=1,2)$.

Proof: By symmetry, it suffices to prove the claim for $l=1$. Let
$\underset{\sim}{a_0},...,\underset{\sim}{a_{n-1}} \in \Gamma_{b_3} \setminus \Gamma_{b_1}$,
$\underset{\sim}{a}=\underset{l<n}{\cup} \underset{\sim}{a_l}$ and
let $\underset{\sim}{g}$ be a $\mathbb{P}_{b_1}$-name of an increasing
sequence from $\omega^{\omega}$, we shall prove that $\Vdash_{\mathbb{P}_{b_3}} "\underset{\sim}{u}:=\{i : \underset{\sim}{a} \cap [g(i),g(i+1))=\emptyset\}$
is infintie$"$. There are $\underset{\sim}{a_l'} \in \Gamma_{a_2} \setminus \Omega_{a_2}$
$(l<n)$ such that $\underset{l<n}{\wedge}f_2(\underset{\sim}{a_l'})=\underset{\sim}{a_l}$,
let $\underset{\sim}{a'}=\underset{l<n}{\cup}\underset{\sim}{a_l'}$.
Let $(\underset{\sim}{m_i} : i<\omega)$ be the $\mathbb{P}_{a_1}$-name
for $f_1^{-1}((\underset{\sim}{g}(i) : i<\omega))$. Let $(p_1,p_2) \in \mathbb{P}_{b_3}$
and $n_*<\omega$, we shall find $(q_1,q_2) \in \mathbb{P}_{b_3}$
above $(p_1,p_2)$ and $n>n_*$ such that $(q_1,q_2) \Vdash_{\mathbb{P}_{b_3}} "n\in \underset{\sim}{u}"$.
We can choose $(p_{1,i},m_{1,i} : i<\omega)$ by induction on $i<\omega$
such that $p_1 \leq p_{1,i} \in \mathbb{P}_{a_1}$, $p_{1,i} \leq p_{1,i+1}$
and $p_{1,i+1} \Vdash_{\mathbb{P}_{a_1}} "\underset{\sim}{m_i}=m_{1,i}"$.
The rest of the proof is as in the previous subclaim. $\square$

\textbf{Claim 8: }For a dense set of $a\in AP$, $\Vdash_{\mathbb{P}_a} "\Gamma_a$
is mad$"$.

\textbf{Proof: }Let $\lambda_0=|\mathbb{P}_a|$ and $\lambda_1=2^{\lambda_0}$.
Let $\mathbb{R}_1=Col(\aleph_0,\lambda_1)$ and $\mathbb P=\mathbb{P}_a \times \mathbb{R}_1 \in H(\kappa)$.

In $V^{\mathbb P}$, $\aleph_1^{V^{\mathbb P}}=\lambda_1^+$ and $\mathbb{P}_a \cup \mathcal{P}(\mathbb{P}_a)$
is countable, so $(\omega^{\omega})^{V^{\mathbb{P}_a}}$ is countable
and $\Gamma:=\{ \underset{\sim}{\tau} : \underset{\sim}{\tau}$ is
a canonical $\mathbb P-$name of a real such that the function listing
$\underset{\sim}{\tau}$ dominates $(\omega^{\omega})^{V^{\mathbb{P}_a}} \}$
is dense in $[\omega]^{\omega}$. By the density of $\Gamma$, we
can find $\Gamma' \subseteq \Gamma$ such that $\Vdash_{\mathbb P}"\Gamma' \cup \Gamma_a$
is mad$"$. Now let $b=(\mathbb P,\Gamma' \cup \Gamma_a)$, then (ignoring
the obvious clauses) we need to prove that $b$ satisfies definition
3(a)(4) and that $a\leq b$ (for which we need to prove that the requirement
from 3(b)(3) is satisfied). We shall prove that $a$ and $b$ satisfy
requirement 3(b)(3), the proof that $b$ satisfies 3(a)(4) is similar.
We shall work in $V^{\mathbb{P}_b}$. Let $\underset{\sim}{a_0},...,\underset{\sim}{a_n} \in \Gamma_b \setminus \Gamma_a$
and let $\underset{\sim}{a}=\underset{l\leq n}{\cup}\underset{\sim}{a_l}$.
Suppose that $(m_i : i<\omega) \in V^{\mathbb{P}_a}$ is increasing,
choose a sequence $(i(k) : k<\omega) \in V^{\mathbb{P}_a}$ such that
$i(k+1)>m_{i(k)+1}+i(k)+(n+1)k$ and let $m_k'=m_{i(k)+1}$ $(k<\omega)$.
For each $l\leq n$, the set $\underset{\sim}{u_l}=\{k<\omega : f_{\underset{\sim}{a_l}}(k)>m_{i(k+1)}\}$
is cofinite (by the definition of $\Gamma$). Therefore, for every
$k$ large enough, $|\underset{\sim}{a_l} \cap m_{i(k+1)}|<k$ (for
every $l\leq n$), hence $|\underset{\sim}{a} \cap m_{i(k+1)}|<(n+1)k$.
For each such $k$, $|\{i : i\in [i(k),i(k+1)) \wedge \underset{\sim}{a} \cap [m_i,m_{i+1}) \neq \emptyset\}|<(n+1)k$.
As $i(k+1)-i(k)>(n+1)k$, there is $i\in [i(k),i(k+1))$ such that
$\underset{\sim}{a} \cap [m_i,m_{i+1})=\emptyset$. Therefore, $f_{\underset{\sim}{a}}$
is not dominated by a real from $V^{\mathbb{P}_a}$. $\square$

\textbf{Claim 9: }For every $a\in AP$ and a $\mathbb{P}_a$-name
$\underset{\sim}{r}$ of a member of $[\omega]^{\omega}$, there is
$b\in AP$ above $a$ such that $\Vdash_{\mathbb{P}_b} "$there is
$\underset{\sim}{s} \in \Gamma_b$ such that $|\underset{\sim}{r} \cap \underset{\sim}{s}|=\aleph_0"$.

\textbf{Proof: }Follows directly from Claim 8. $\square$

\textbf{Observation 10: }Let $\mathbb{Q}$ be a forcing notion from
$\bold Q$.\textbf{ }Assume that $a_0 \leq a_l$, $\underset{\sim}{\eta_l}$
is a $\mathbb{P}_{a_l}$-name of a $\mathbb Q-$generic real over
$V^{\mathbb{P}_{a_0}}$ $(l=1,2)$, and $\mathbb{P}_{a_0} \star \underset{\sim}{\eta_1}$
is isomorphic to $\mathbb{P}_{a_0} \star \underset{\sim}{\eta_2}$
over $\mathbb{P}_{a_0}$ (so wlog they're equal to each other and
we may denote the generic real by $\underset{\sim}{\eta}$). By Claim
6, there is $a_0' \in AP$ such that $a_0 \leq a_0' \leq a_l$ $(l=1,2)$,
$\mathbb{P}_{a_0'}=\mathbb{P}_{a_0} \star \underset{\sim}{\eta}$
and $\Gamma_{a_0'}=\Gamma_{a_0}$. By Claim 7, there are $b_l$ $(l\leq 3)$
and $g_l$ $(l\leq 2)$ as there for $(a_0',a_1,a_2)$ here. $\square$

\textbf{Definition 11: }Let $H\subseteq AP$ be generic over $V$
and let $V_1=V[H]$. In $V_1$, let $\underset{\sim}{\mathbb P}[H]$
be $\underset{a\in H}{\cup}\mathbb{P}_a$.

\textbf{Claim 12: }$\Vdash_{AP} "\underset{\sim}{\mathbb P} \models \kappa-cc"$.

\textbf{Proof: }Suppose towards contradiction that $\Vdash_{AP} "\underset{\sim}{I} \subseteq \underset{\sim}{\mathbb P}$
is a maximal antichain of cardinality $\kappa"$. Choose by induction
on $\alpha<\kappa$ a sequence $(a_{\alpha},p_{\alpha} : \alpha<\kappa)$
such that:

a. $a_{\alpha} \in AP$.

b. $(a_{\beta} : \beta<\alpha)$ is $\leq_{AP}$-increasing cotinuous.

c. $a_{\beta+1} \Vdash_{AP} "p_{\beta} \in \underset{\sim}{I} \setminus \{p_{\gamma} : \gamma<\beta\}"$.

d. $p_{\beta} \in \mathbb{P}_{\beta+1}$.

For every $\alpha<\kappa$, there is $q_{\alpha} \in \mathbb{P}_{a_{<\alpha}}:=\underset{\gamma<\alpha}{\cup}\mathbb{P}_{a_{\gamma}}$
such that $p_{\alpha}$ is compatible with every $r\in \mathbb{P}_{a_{<\alpha}}$
above $q_{\alpha}$. Let $\gamma(\alpha)<\alpha$ be the least $\gamma$
such that $q_{\alpha} \in \mathbb{P}_{a_{\gamma}}$. For some $\gamma(*)<\kappa$,
$S:=\{ \alpha: \gamma(\alpha)=\gamma(*)\}$ is stationary. As $|\mathbb{P}_{a_{\gamma(*)}}|<\kappa$,
there is $S'\subseteq S$ of cardinality $\kappa$ such that $\alpha_1<\alpha_2 \in S' \rightarrow q_{\alpha_1}=q_{\alpha_2}$,
which leads to a contradiction. $\square$

\textbf{Definition 13: }Let $V_1$ be as in Definition 11 and let
$G\subseteq \underset{\sim}{\mathbb{P}[H]}$ be generic over $V_1$,
we shall denote $V[H,G]$ by $V_2$.

\textbf{Caim 14: }Every real in $V_2$ is from $V_1[G\cap \mathbb{P}_a]$
for some $a\in H$.

\textbf{Proof: }Let $\underset{\sim}{r}$ be a $AP \star \underset{\sim}{\mathbb P}$-name
of a real. By Claim 12, $\underset{\sim}{\mathbb P}[H] \models \kappa-cc$
in $V_1$. Therefore, for every $n<\omega$ there are $AP-$names
$\bar{p_n}=(\underset{\sim}{p_{n,\alpha}} : \alpha<\underset{\sim}{\alpha_n})$
and $\bar{t_n}=(\underset{\sim}{t_{n,\alpha}} : \alpha<\underset{\sim}{\alpha_n})$
such that:

a. $\underset{\sim}{\alpha_n}<\kappa$.

b. $\bar{p_n}$ is a maximal antichain in $\underset{\sim}{\mathbb{P}}[H]$.

c. $\underset{\sim}{t_{n,\alpha}}$ is a $\underset{\sim}{\mathbb{P}}[H]-$name
of an element of $\{0,1\}$.

d. $\underset{\sim}{p_{n,\alpha}} \Vdash "n\in \underset{\sim}{r}$
iff $\underset{\sim}{t_{n,\alpha}}=1"$.

For every $n<\omega$ and $\alpha<\underset{\sim}{\alpha_n}$, there
is $\underset{\sim}{a_{n,\alpha}} \in \underset{\sim}{H}$ such that
$\underset{\sim}{p_{n,\alpha}} \in \mathbb{P}_{\underset{\sim}{a_{n,\alpha}}}$.
Now let $a_0 \in AP$, we can find $\leq_{AP}$-increasing sequence
$(a_n : n<\omega)$ such that $a_{n+1} \Vdash "\underset{\sim}{\alpha_n}=\alpha_n^*"$
for some $\alpha_n^*<\kappa$. Let $a_{\omega} \in AP$ be an upper
bound, and now choose an increasing sequence $(a_{\omega+\alpha} : \alpha \leq \underset{n<\omega}{\Sigma} \alpha_n^*)$
by induction on $\alpha \leq \underset{n<\omega}{\Sigma} \alpha_n^*$
such that for every $n<\omega$ and $\beta<\alpha_n^*$, $a_{\omega+\underset{l<n}{\Sigma} \alpha_l^* +\beta+1} \Vdash "\underset{\sim}{a_{n,\beta}}=a_{n,\beta}^*$
and $\underset{\sim}{p_{n,\beta}}=p_{n,\beta}^*"$. We may assume
wlog that $a_{n,\beta}^* \leq_{AP} a_{\omega+\underset{l<n}{\Sigma} \alpha_l^* +\beta+1}$,
so $p_{n,\beta}^* \in \mathbb{P}_{a_{\omega++\underset{l<n}{\Sigma} \alpha_l^* +\beta+1}}$.
It's now easy to see that $\underset{\sim}{r}$ is a $\mathbb{P}_{a_{\omega+\underset{n<\omega}{\Sigma} \alpha_n^*}}$-name.
$\square$ 

\textbf{Theorem 15: }a. In $V_2$, let $\mathcal A=\{\underset{\sim}{a}[G] : \underset{\sim}{a} \in \Gamma_b$
for some $b\in H\}$ and let $V_3=HOD(\mathbb R, \mathcal A)$, then
$V_3 \models ZF+DC+"$there exists a mad family$"+"$all sets of reals
are $\mathbb{Q}-$measurable for every $\mathbb Q \in \bold Q"$.

b. $ZF+DC+"$every set of reals is Lebesgue measurable$"+"$there
exists a mad family$"$ is consistent relative to an inaccessible
cardinal.

\textbf{Proof: }a. The existence of a mad family follows by Claim
8. $\mathbb Q-$measurability for $\mathbb Q \in \bold Q$ follows
from Claim 14 and Observation 10 as in Solovay's proof.

b. Apply the previous clause to $\mathbb Q=$Random real forcing.
$\square$

As a corollary to the above theorem, we obtain an answer to a question
of Henle, Mathias and Woodin from {[}HMW{]}: 

\textbf{Corollary 16 $(ZF+DC)$: }The existence of a mad family does
not imply that $\aleph_1 \leq \mathbb R$.

\textbf{Proof: }By Theorem 15 (applied to Random real forcing) and
the fact that the existence of an $\omega_1$-sequence of distinct
reals implies the existence of a non-Lebesgue measurable set of reals
(see {[}Sh176{]}). $\square$

\textbf{Remark}: The above result was also obtained by Larson and
Zapletal in {[}LZ{]} assuming the existence of a proper class of Woodin
cardinals.

We conclude with a somewhat surprising observation, showing that the
analog of Theorem 15 fails at the lower levels of the projective hierarchy:

\textbf{Observation 17: }If every $\Sigma^1_3$ set of reals is Lebesgue
measurable, then there are no $\Sigma^1_2$-mad families.

\textbf{Proof: }By $[Sh176]$, $\Sigma^1_3$-Lebesgue measurability
implies that $\omega_1^{L[x]}<\omega_1$ for every $x\in \omega^{\omega}$.
By Theorem 1.3(2) in {[}To{]}, it follows that there are no $\Sigma^1_2$-mad
families. $\square$

\noindent \begin{center}
\textbf{\large On a question of Enayat}
\par\end{center}{\large \par}

We now address a question asked by Ali Enayat in {[}En{]}. The question
is motivated by the problem of understanding the relationship between
Freiling's axiom of symmetry, the continuum hypothesis and the Lebesgue
measurability of all sets of reals (see discussion in {[}Ch{]}). 

As with the previous results, we were informed by Paul Larson that
the following results can also be obtained under the assumption of
a proper class of Woodin cardinals using the arguments from {[}LZ{]}.

\textbf{Definition 18: }a. Let $WCH$ (weak continuum hypothesis)
be the statement that every uncountable set of reals can be put into
1-1 correspondence with $\mathbb R$.

b. Let $AX$ (Freiling's axiom of symmetry) be the following statement:
Let $\mathcal F$ be the set of functions $f: [0,1] \rightarrow \mathcal{P}_{\omega_1}([0,1])$,
then for every $f\in \mathcal F$ there exist $x,y \in [0,1]$ such
that $x\notin f(y)$ and $y\notin f(x)$.

Remark: The term $WCH$ has a different meaning in several papers
by other authors.

\textbf{Theorem 19: }$ZF+DC+\neg WCH+"$every set of reals is Lebesgue
measurable$"$ is consistent relative to an inaccessible cardinal.

\textbf{Proof: }Let $V_3$ be the model from Theorem 15(b), we shall
prove that $V_3 \models \neg WCH$ by showing that there is no injection
from $\mathbb R$ to the mad family $\mathcal A$. Suppose toward
contradiction that for some $(a,\underset{\sim}{p}) \in AP \star \underset{\sim}{\mathbb P}$
(where $\underset{\sim}{\mathbb P}$ is as in Definition 11), a canonical
name for a real $\underset{\sim}{r}$ and a first order formula $\phi(x,y,z,\mathcal A)$,
$(a,\underset{\sim}{p})\Vdash "\phi(x,y,\underset{\sim}{r},\mathcal A)$
defines an injection $F_{\underset{\sim}{r}}$ from $\mathbb R$ to
$\mathcal A"$. We may assume wlog that $\underset{\sim}{r}$ is a
canonical $\mathbb{P}_a$-name. We may also assume wlog that, for
every $\underset{\sim}{s} \in \Gamma_a$, $(a,\underset{\sim}{p})\Vdash "$if
$\underset{\sim}{s} \in Ran(F_{\underset{\sim}{r}})$, then $\underset{\sim}{s}=F_{\underset{\sim}{r}}(t)$
for some $t\in \mathbb{R}^{V^{\mathbb{P}_a}}"$. This is possible
as $|\Gamma_a|<\kappa$, so we may construct an increasing sequence
$(a_{\gamma} : \gamma<\beta)$ of length $<\kappa$, such that $a_0=a$
and such that the upper bound $(a_{\beta},\Gamma_{a_\beta})$ satisfies
the above requirement. $((a_{\beta},\Gamma_a),\underset{\sim}{p})$
is then as required. By increasing $a$, we may assume wlog that $\underset{\sim}{p}$
is an object $p$ (and not just an $AP$-name) from $\mathbb{P}_a$.
Now let $a_2 \in AP$ be defined as $a_2=(\mathbb{P}_a \star Cohen, \Gamma_a)$
and let $\underset{\sim}{\eta}$ be the $\mathbb{P}_{a_2}$-name for
the Cohen real. There are $a_3 \in AP$ and a name $\underset{\sim}{\nu}$
such that $a_2 \leq a_3$ and $a_3 \Vdash "p \Vdash "\phi(\underset{\sim}{\eta},\underset{\sim}{\nu},\underset{\sim}{r},\mathcal A)""$,
so $\underset{\sim}{\nu} \in \mathcal A$, and by the injectivity
of $F_{\underset{\sim}{r}}$, $\underset{\sim}{\nu} \notin \Gamma_a$.
We may assume wlog that $\underset{\sim}{\nu} \in \Gamma_{a_3}$.

Let $a_4$ be the amalgamation of two copies of $a_3$ over $a_2$
(i.e. as in the proof of Claim 7) and let $f_0: \mathbb{P}_{a_3} \rightarrow \mathbb{P}_{a_4}$
and $f_1: \mathbb{P}_{a_3} \rightarrow \mathbb{P}_{a_4}$ be the corresponding
complete embeddings. As the amalgamation is over $a_2$, it follows
that $f_0(\underset{\sim}{\eta})=f_1(\underset{\sim}{\eta})$ and
$f_0(\underset{\sim}{r})=f_1(\underset{\sim}{r})$, and by the argument
from the proof of Claim 7 (Subclaim 5), $f_0(\underset{\sim}{\nu}) \neq f_1(\underset{\sim}{\nu})$.
As $f_l$ $(l=0,1)$ are isormorphisms between $a_3$ and $f_l(a_3) \leq a_4$
such that $f_l \restriction \mathbb{P}_{a_2}=Id$, they induce an
automorphism of $(AP,\leq_{AP})$ mapping $a_3$ to $f_l(a_3)$ and
$a_2$ to itself. Therefore, $a_4 \Vdash "" f_0(p) \Vdash"\phi(f_0(\underset{\sim}{\eta}),f_0(\underset{\sim}{\nu}),f_0(\underset{\sim}{r}),\mathcal A)""$, 

$a_4 \Vdash ""f_1(p) \Vdash "\phi(f_1(\underset{\sim}{\eta}),f_1(\underset{\sim}{\nu}),f_1(\underset{\sim}{r}),\mathcal A)""$
and $f_0(p)=f_1(p)$, a contradiction. $\square$

\textbf{Theorem 20: }$WCH$ is independent of $ZF+DC+AX+"$all sets
of reals are Lebesgue measurable$"$.

\textbf{Proof: }By {[}We{]}, $AX$ is implied by $ZF+DC+"$all sets
of reals are Lebesgue measurable$"$. Therefore, $AX$ holds in the
model $V_3$ from Theorem 15(b) and in Solovay's model. By Corollary
19, $V_3 \models \neg WCH$. By the fact that all sets of reals in
Solovay's model have the perfect set property, it follows that $WCH$
holds in that model. $\square$

\textbf{\large References}{\large \par}

{[}Ch{]} Timothy Chow, Question about Freiling's axiom of symmetry,
cs.nyu.edu/pipermail/fom/2011-August/015676.html

{[}En{]} Ali Enayat, Lebesgue measurability and weak CH, mathoverflow.net/questions/72047/lebesgue-measurability-and-weak-ch

{[}HwSh1090{]} Haim Horowitz and Saharon Shelah, Can you take Toernquist
inaccessible away?, arXiv:1605.02419 

{[}HMW{]} James Henle, Adrian R. D. Mathias and W. Hugh Woodin, A
barren extension. Methods in mathematical logic, Lecture Notes in
Mathematics 1130, pages 195-207. Springer Verlag, New York, 1985

{[}LZ{]} Paul Larson and Jindrich Zapletal, Canonical models for fragments
of the axiom of choice, users.miamioh.edu/larsonpb/lru5.pdf

{[}Ma{]} A. R. D Mathias, Happy families, Ann. Math. Logic \textbf{12
}(1977), no. 1, 59-111

{[}NN{]} Itay Neeman and Zach Norwood, Happy and mad families in $L(\mathbb R)$,
math.ucla.edu/\textasciitilde{}ineeman/hmlr.pdf

{[}Sh176{]} Saharon Shelah, Can you take Solovay's inaccessible away?
Israel J. Math. 48 (1984), no. 1, 1-47

{[}Sh218{]} Saharon Shelah, On measure and category, Israel J. Math.
52 (1985) 110-114

{[}So{]} Robert M. Solovay, A model of set theory in which every set
of reals is Lebesgue measurable, AM 92 (1970), 1-56

{[}To{]} Asger Toernquist, Definability and almost disjoint families,
arXiv:1503.07577

{[}We{]} Galen Weitkamp, The $\Sigma^1_2$-theory of axioms of symmetry,
J. Symbolic Logic 54 (1989), no. 3, 727-734

$\\$

(Haim Horowitz) Einstein Institute of Mathematics

Edmond J. Safra campus, 

The Hebrew University of Jerusalem.

Givat Ram, Jerusalem, 91904, Israel.

E-mail address: haim.horowitz@mail.huji.ac.il

$\\$

(Saharon Shelah) Einstein Institute of Mathematics

Edmond J. Safra campus, 

The Hebrew University of Jerusalem.

Givat Ram, Jerusalem, 91904, Israel.

Department of Mathematics

Hill Center - Busch Campus,

Rutgers, The State University of New Jersey.

110 Frelinghuysen road, Piscataway, NJ 08854-8019 USA

E-mail address: shelah@math.huji.ac.il
\end{document}